\documentclass[12pt]{article}
\usepackage[utf8]{inputenc}
\usepackage{graphicx}
\usepackage[T2A]{fontenc}
\usepackage[english,russian]{babel}
\usepackage{amsmath,amssymb,euscript,amsthm,amsfonts,mathrsfs,amscd,amsthm}
\usepackage{color}
\usepackage{floatflt}
\usepackage{mathrsfs} 
\usepackage[ left=35mm, top=25mm, right=35mm, bottom=50mm]{geometry}
\usepackage{amsthm}
\newtheorem{theorem}{Теорема}
\newtheorem*{myth}{Теорема}
\newtheorem*{mylem}{Лемма}

\newtheorem*{mydef}{Определение}

\newtheorem*{myrem}{Замечание}

\newcommand{\bk}{\color{black}}
\newcommand{\gr}{\color{black}}
\newcommand{\bb}{\color{black}}
\newcommand{\rr}{\color{black}}
\newcommand{\dd}{\ensuremath{\displaystyle}}

\newcommand{\supl}{\ensuremath{\sup\limits}}
\newcommand\intl{\ensuremath{\int\limits}}
\newcommand{\suml}{\ensuremath{\sum\limits}}
\newcommand{\bD}{\stackrel{\mathscr{D}}{=}}
\newcommand{\1}{\ensuremath{\mathbf{1}}}
\newcommand{\PPP}{\mathscr{P}}
\newcommand{\FFF}{\mathscr{F}}

\newcommand{\XX}{\mathscr{X}}

\newcommand{\SSS}{\mathscr{S}}

\newcommand{\UU}{\mathscr{U}}
\newcommand{\EEE}{\mathscr{E}}

\newcommand{\ud}{\,\mathrm{d}}

\newcommand{\bd}{\stackrel{\text{\rm def}}{=\!\!\!=}}
\newcommand{\PP}{\ensuremath{\mathbf{P}}}

\newcommand{\EE}{\ensuremath{{\mathbb{E}}}}
\newcommand{\const}{\ensuremath{\mbox{const }}}

\newcommand*{\QEDB}{\hfill \ensuremath{\square}}%

\newcommand*{\TR}{\hfill \ensuremath{\triangleright}}

\title{Об оценке скорости сходимости для регенерирующих процессов в ТМО и смежных задачах. }
\author{М.П.~Фархадов, Г.А.~Зверкина\footnote{Работа автора поддержана грантом РФФИ № 20-01-00575А}
\\
\small Институт проблем управления им. В. А. Трапезникова РАН\\
\small 
(Institute of Control Sciences V.A. Trapeznikov, Russian Academy of Sciences)}

\begin{document}
\maketitle
\begin{abstract}
   Цель этой статьи – показать, как метод склеивания может применяться для получения строгих оценок сверху для скорости сходимости распределения регенерирующих марковских процессов к стационарному распределению в случаях, когда рассматриваемый регенерирующий марковский процесс эргодичен. 
   Этот метод может применяться к СМО, СеМО и системам надёжности. 
   Статья представляет собой изложение пленарного доклада авторов на VIII Международной молодежной научной конференции «Математическое и программное обеспечение информационных, технических и экономических систем» (г. Томск, 26-30 мая 2021 года).
\end{abstract}
{\bf Ключевые слова:} {\it 
теория массового обслуживания; регенерирующие марковские процессы; метод склеивания для кусочно-линейных марковских процессов; метрика полной вариации.
}
\section*{Введение}

В теории массового обслуживания (ТМО) и в смежных задачах очень важно знать числовые характеристики рассматриваемой системы – как в стационарном, так и в достационарном режиме.
В ряде случаев такие характеристики могут быть вычислены, но это возможно для ограниченного количества реальных моделей.
Однако в большинстве случаев возможно вычисление или оценка стационарных значений характеристик исследуемых моделей.
Поведение подавляющего количества систем массового обслуживания (СМО), систем надёжности и сетей массового обслуживания (СеМО) может быть описано с помощью линейчатых марковских процессов, которые во многих случаях являются регенерирующими.
В том случае, когда период регенерации марковского процесса имеет конечное среднее значение, этот процесс эргодичен.
Метод склеивания может применяться для получения строгих оценок сверху для скорости сходимости распределения регенерирующих марковских процессов (РМП) к стационарному распределению (естественно, в тех случаях, когда рассматриваемый РМП эргодичен). Этот метод может применяться к анализу систем массового обслуживания, сетей массового обслуживания и систем надёжности.

К сожалению, в отечественной литературе практически нет информации о методе склеивания (см., например, \cite{1.}, \cite{5.}), хотя отечественные учёные использовали этот метод (\cite{6.} и др.).

\section{Кусочно-линейные Марковские процессы \linebreak(КЛМП)}

Поведение систем массового обслуживания (СМО), сетей массового обслуживания (СеМО) и систем надёжности определяется, как правило, последовательностями случайных величин (сл.в.), представляющих собой времена работы или обслуживания требований (восстанавливаемых элементов или узлов СеМО), интервалы между поступлениями требований в исследуемую систему, случайные времена ожидания обслуживания (в случае “нетерпеливых” требований) и пр.

Такого сорта поведение технических систем принято описывать кусочно-линейными (или линейчатыми) процессами, введёнными в \cite{4.}.

Следуя \cite{4.}, определим кусочно-линейный марковский процесс следующим образом.

Кусочно-линейным Марковским процессом называется случайный процесс $X_t = \left(\nu(t), \vec\upsilon(t)\right)$, определяемый следующим образом.
Пространство состояний процесса $\XX$ -- это множество пар $\left(i,\vec \xi_i\right)$, где $i$ -- элемент конечного или счётного пространства, а $\vec\upsilon_i$ -- вектор $=(\upsilon_1, \ldots \upsilon_{|i|})$, где $|i|\geqslant 0$ является ``размером'' {\it базового состояния} $i$; $\upsilon_j\geqslant O$. 
Поведение процесса $X_t$ описывается так. 

Пусть $X_t = (i,\vec y),$ $\vec y = (y_l, \ldots y_{|i|})$. 
Тогда с вероятностью $\lambda_{i}(X_t)\times p_{ij}(X_t) \ud t$
за время $\ud t$ произойдёт случайный переход $X_t$ в {\it базовое состояние} $j$.
После перехода в состояние $j$ новое значение компоненты $ \vec\upsilon(t)$ случайно и определяется измеримой по $\vec y$ функцией распределения (ф.р.)
$$
B_{ij}^{(0)}(\vec x| \vec y) =\PP\{ \vec\upsilon(t + \ud t) < \vec x | \vec\upsilon(t) = \vec y,\, \nu(t ) = i, \,\nu(t + \ud t) = j\}
$$

Вероятность того, что за малое время $h$ произойдёт более одного случайного перехода есть $o(h)$. 
При этом, если $\nu(t + \ud t) = \nu(t)$, то $\vec\upsilon(t + \ud t) =\vec\upsilon(t) + \vec \alpha_i \ud t$, где $\vec \alpha_i= (\alpha_{il} , \ldots , \alpha_{i|i|})$ -- вектор с неотрицательными компонентами, если на интервале $(t, t + \ud t)$ не было случайных переходов (см. \cite{4.}).

Для СМО номер $i$ может обозначать количество находящихся в СМО требований, а вектор $\vec\upsilon$ состоит из: времени, прошедшего со времени прихода последнего требования, а также из времён пребывания имеющихся требований с СМО (или стоимость уже проведённых работ по их обслуживанию и проч.).

В СеМО $i$ может нумеровать все возможные состояния СеМО (количество требований в каждом узле и в очередях на этих узлах), а вектор $\vec\upsilon$ может характеризовать состояние (прошедшее время обслуживания или ожидания обслуживания) находящихся на обслуживании на каждом узле требований.

Таких вариантов описания моделей ТМО и смежных задач множество.

Отметим, что в случае, когда $\lambda(X_t)\equiv\lambda(\nu_t)$ и $p_{ij}(X_t)\equiv p(i,j)$, процесс $X_t$ представляет собой цепь Маркова в непрерывном времени (см. \cite{4.}, \cite{7.}).

\subsection{\it Вложенная цепь Маркова}

КЛМП имеет вложенную {\it неоднородную} цепь Маркова -- с вероятностями перехода $p_{ij}(X(t))$ (обычно рассматривается случай, когда $p_{ij}(X(t))\equiv p(i,j)=\const\!\!$).

Мы предполагаем, что эта Марковская цепь неприводима и положительно возвратна. 
При этом $\lambda_{i}(X_t)$ -- это интенсивность окончания времени пребывания в текущем состоянии $\nu(t)=i$.
Как уже говорилось, в случае, если $\lambda_{i}(X_t)\equiv\lambda (\nu(t))$, КЛМП -- это цепь Маркова в непрерывном времени.

Хорошо известно, что поведение цепи Маркова в непрерывном времени описывается уравнениями Колмогорова.

Если среднее время пребывания во всех состояниях цепи Маркова конечно, что соответствующий КЛМП эргодичен(см. \cite{7.}).

\section{Регенерирующие процессы}

Поведение большинства изучаемых СМО, СеМО и систем надёжности описывается регенерирующими процессами. 

\begin{mydef}
[Регенерирующего процесса]
Процесс $(X_t,\,t\geqslant 0)$, заданный на вероятностном пространстве $(\Omega,\mathscr{F} , \PP )$, с измеримым пространством состояний $(\mathscr{X} , \mathscr{B} (\mathscr{X} ))$ называется регенерирующим, если существует возрастающая последовательность $\left\{ \theta_n\right\} $ $(n\in\mathbb{Z}_+)$ Марковских моментов по отношению к фильтрации $\mathscr{F} _{t\leqslant 0}$ таких, что последовательность 
$$
\left\{ \Theta_n\right\} =\left\{ B_{t+\theta_{n-1}}-B_{\theta_{n-1}}, \theta_n- \theta_{n-1},t\in[\theta_{n-1},\theta_n)\right\} ,\quad n\in\mathbb{N}
$$
состоит из независимых одинаково распределённых (н.о.р.) случайных элементов, заданных на $(\Omega,\mathscr{F} , \PP )$. 
Если $\theta_0\neq 0$, то регенерирующий процесс $(X_t,\,t\geqslant 0)$ называется процессом с запаздыванием.
Моменты $\left\{ \theta_n\right\} $ называются моментами (точками) регенерации.
{\hfill \ensuremath{\triangleright}}
\end{mydef}
Обозначим $\xi_n\stackrel {\text{\rm def}}{=\!\!\!=} \theta_n- \theta_{n-1}$, и пусть $F(s)=\PP \left\{ \xi_n\leqslant s\right\} $ -- ф.р. длины периода регенерации; везде далее мы предполагаем, что {\it распределение $F$ нерешётчато} и, более того, для простоты будем считать, почти всюду при $t>0$ существует положительная плотность распределения.

В ТМО моментами регенерации обычно являются моменты начала периода занятости (при экспоненциальном входящем потоке это могут быть и моменты окончания периодов занятости).
В теории надёжности при анализе поведения {\it одного } восстанавливаемого элемента точки регенерации -- это моменты отказов или восстановлений.
При анализе СеМО с экспоненциальными распределениями времени пребывания требований в узлах СеМО и входящих экспоненциальных потоках КЛМП представляет собой цепь Маркова в непрерывном времени; в этом случае в качестве точки регенерации можно выбрать любое {\it возвратное} базовое состояние. 

\subsection {\it Наша цель}

\quad Пусть $\PPP_t$ -- распределение КЛМП $(X_t,\,t\geqslant 0)$ в момент $t$.

\quad Если $\EE\,\xi_n<\infty$, то существует и единственно стационарное инвариантное распределение $\PPP$ такое, что $\PPP_t\Longrightarrow \PPP$.

Наша цель -- вычисление строгой оценки сверху для скорости сходимости $\PPP_t\Longrightarrow \PPP$.

Для этого будет использован {\it метод склеивания} Марковских процессов, т.е. исследуемые регенерирующие процессы $(X_t,\,t\geqslant 0)$ должны быть Марковскими.
Если процесс $(X_t,\,t\geqslant 0)$ не является Марковским, то можно расширить пространство его состояний $\mathscr{X}$ таким образом, что на новом (расширенном) пространстве состояний $\overline{\mathscr{X}}$ ``уточнённый'' процесс $(\overline X_t,\,t\geqslant 0)$ окажется марковским.

Например, в состояние регенерирующего процесса $X_t$, $t\in[\theta_{n-1},\theta_n)$ можно включить полную историю этого процесса на интервале $[\theta_{n-1},t]$: процесс $\overline X_t\stackrel {\text{\rm def}}{=\!\!\!=} \left\{ X_s, s\in[\theta_{n-1},t] |t<\theta_n\right\} $ -- марковский и регенерирующий на расширенном пространстве состояний $\overline{ \mathscr{X}} $.

\subsection{\it Что уже известно о регенерирующих процессах (см. \cite{2.})}

\begin{myth}
Если для некоторого $\alpha>0$ выполнено условие 
$$
\mathbf{E}\,e^{\alpha \xi}<\infty,
$$
то для любого $\beta<\alpha$ существует постоянная $C(\beta)$ такая, что для всех $A\in\sigma(\mathscr{X})$ и всех $t>0$ верно:
$
|\PPP_t( A)-\PPP(A) |<C(\beta)\exp(-\beta t).
$ \QEDB
\end{myth}
Оценки для постоянной $C(\beta)$ неизвестны.

\begin{myth}
Если для некоторого $k>1$ выполнено условие 
$$
\qquad\mathbf{E}\,\xi^k<\infty,
$$
то для любого $\kappa\leqslant k-1$ существует постоянная $C(\kappa)$ такая, что\linebreak
$
|\PPP_t( A)-\PPP(A) |<C(\kappa)t^{-\kappa}.
$ \QEDB
\end{myth}
Оценки для постоянной $C(k)$ неизвестны.

\subsection{\it Вложенный процесс восстановления}

Рассмотрим регенерирующий процесс $X_t$.
Этот процесс имеет вложенный процесс восстановления $N_t$

\begin{figure}[h]
\center{\includegraphics[width=400pt]{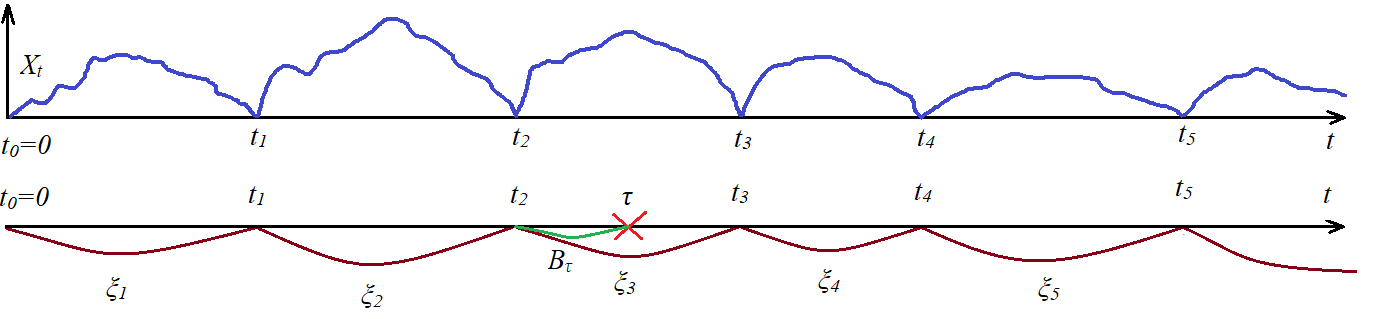}}
\caption{Регенерирующий процесс ${{X}_{t}}$ и вложенный процесс восстановления.}
\label{fig0}
\end{figure}
Распределение процесса $X_\tau$ -- это функция от значения случайной величины $B_\tau$ -- обратного времени восстановления.

Значит, если известны оценки скорости сходимости распределения обратного времени восстановления процесса восстановления $ N_t $, то известны оценки скорости сходимости распределение регенерирующего процесса $ X_t $.

\subsection{\it Процесс восстановления и неравенство Лордена}

Рассмотрим процесс восстановления $N _t\bd\dd\suml _{i=1} ^\infty \1\left \{ \suml _{s=1} ^i \xi _k\leqslant t\right \} $, где
$\left \{\xi _1, \xi _2, ...\right \} $ -- н.о.р. положительные случайные величины (сл.в.) с ф.р. $F(s)$.
$N _t$ -- это считающий процесс, который меняет своё значение в моменты времени $t _k=S _k\bd\dd\suml _{j=1} ^s \xi _j$.
Эти моменты $t _k$ называются моментами восстановления.

\begin{figure}[h]
\centering
\begin{picture}(300,45)
\put(0,20){\vector(1,0){300}}
\put(0,18){\line(0,1){4}}
\thicklines
{\qbezier(0,20)(30,40)(60,20)}
{\qbezier(60,20)(80,40)(100,20)}
\qbezier(100,20)(150,40)(200,20)
\qbezier(200,20)(240,40)(280,20)
\put(280,30){\ldots \ldots \ldots}
\put(58,10){$t _1$}
\put(98,10){$t _2$}
\put(198,10){$t _3$}
\put(278,10){$t _4$}
\put(27,40){$\xi _1$}
\put(77,40){$\xi _2\bD \xi _1$}
\put(147,40){$\xi _3\bD \xi _1$}
\put(237,40){$\xi _4\bD \xi _1$}
\put(-5,10){$t _0=0$}
\put(148,10){$t$}
\put(195,15){\line(1,1){10}}
\put(195,25){\line(1,-1){10}}
\rr
\qbezier(150,20)(175,0)(200,20)
\qbezier(150,20)(175,1)(200,20)
\qbezier(150,20)(175,2)(200,20)
\put(170,0){$W _t$}
\bb
\qbezier(100,20)(125,0)(150,20)
\qbezier(100,20)(125,1)(150,20)
\qbezier(100,20)(125,2)(150,20)
\put(120,0){$B _t$}
%
\end{picture}
\caption{$B _t$ -- перескок, $W _t$ -- недоскок в момент времени $t$ (см. \cite{4.}, \cite{8.})}
\label{fig1}
\end{figure} 
На Рис. \ref{fig1} показаны обратное время восстановления или перескок $B _t$ и прямое время восстановления или недоскок $W _t$ в {\bf фиксированный} момент времени $t$: 
$B _t\bd t-S _{N _t}; \quad W _t\bd S _{N _t+1}-t.$

Напомним, что в наших условиях предельное распределение величины перескока и недоскока при $t\to \infty $ таково: $$\tilde{F}(s)=1-\frac{\int\limits_{0}^{s}{(1-F(}u))\text{d}u}{\int\limits_{0}^{\infty }{(1-F(}u))\text{d}u}=1-\frac{\int\limits_{0}^{s}{(1-F(}u))\text{d}u}{\mathbb{E}\xi }$$ –- см.
\cite{8.}.

\begin{myth}[Г.Лордена, см. \cite{11.}]
Неравенство Лордена даёт оценку сверху математического ожидания (м.о.) перескока:
$$
\EE\, B _t \leqslant \frac{\EE\, \xi ^2}{\EE\, \xi}\bd \Xi\big(=\mbox{функционал от }F(s)\big).\eqno{\QEDB}
$$
\end{myth}
Это неравенство будет использовано для получения оценки сверху для скорости сходимости распределения перескока к стационарному распределению в метрике полной вариации.

Напомним, что 
\begin{multline*}
\PP\{\mbox{\it имеется хотя бы одно восстановление на интервале }[t, t+\Delta]|{\bb B_t}=y\} =\\
= \frac{F(y+\Delta)- F(y)}{1-F(y)}
=\intl_y^{y+\Delta}\lambda(s)\ud s=\lambda(y)\Delta + o(\Delta), 
\end{multline*}
где 
\begin{equation}\label{1}
\lambda(t)= \dd\frac{F'(t)}{1-F(t)}={\frac{\ud}{\ud t}\Big(-\ln(1-F(t)\Big)},    
\end{equation}
и $\lambda(t)$ называется {\it интенсивностью процесса восстановления}.

При этом 
\begin{equation}\label{2}
{F(s)=1-\exp\left(\intl_0^s(-\lambda(u))\ud u\right)}, \qquad F'(s)=\lambda(s)\exp\left(\intl_0^s(-\lambda(u))\ud u\right).    
\end{equation}

\begin{myrem}
Обычно предполагается, что распределение $F(s)$ абсолютно \linebreak непрерывно. Однако в некоторых приложениях теории восстановления ф.р. сл.в. $\xi_i$ не абсолютно непрерывны.

Т.е. ф.р. сл.в. $\xi_i$ может иметь скачки (естественно, в прикладных задачах теории вероятностей мы не рассматриваем сингулярные сл.в.).

Положим
$
f(s)=
\begin{cases}
F'(s), & \mbox{если } F'(s) \mbox{ существует; }\\
0, & \mbox{в противном случае}.
\end{cases}
$

В дальнейшем полагаем $\lambda(s)\bd \dd\frac{f(s)}{1-F(s)}-\suml_{i}\delta(s-a_i)\ln(F(a_i+0)-F(a_i-0))$, где $\{a_i\}$ -- множество точек разрыва $F(s)$. 
Такое определение интенсивности оставляет в силе первую формулу в (\ref{2})  -- это следует из (\ref{1}). \TR
\end{myrem}

 Заметим, что процесс $\bb B_t$ -- это простейший Марковский КЛМП или {\it линейчатый процесс} (см. \cite{4.}).

\section{Метод склеивания (Марковских процессов) -- см., например, \cite{9.}}

Применение метода склеивания основано на следующем соображении.

{\it Если два однородных независимых Марковских процесса с одними и теми же переходными вероятностями, но с разными начальными условиями, совпадут в некоторый момент $\tau$, то после этого момента $\tau$ их распределения совпадают.
Момент $\tau$ называется временем (моментом) склеивания или склейкой.}

Напомним, что процесс $X_t$ с пространством состояний $\XX$ называется Марковским, если для него выполнено {\it Марковское свойство}:
\\
$\forall A_0,A_1,\ldots,A_n\in\sigma(\XX)$,\; $\forall t_1,t_2,\ldots,t_n, t:$\; $0\leqslant t_i<t_{i+1}<t$ верно равенство:
$$\PP\{X_t\in A_0|X_{t_1}\in A_1, X_{t_2}\in A_2, \ldots, X_{t_n}\in A_n\}=\PP\{X_t\in A_0|X_{t_n}\in A_n\},$$
т.е. поведение процесса после момента $t_n$ зависит только от его состояния в момент $t_n$ и не зависит от того, что было раньше.

Пусть теперь два независимых Марковских процесса $X_t$ и $X_t'$ имеют одинаковые переходные вероятности: $\PP\{X_{t+\theta}\in A|X_t\in B\}\equiv \PP\{X'_{s+\theta}\in A|X_s'\in B\}$ для всех $\theta>0$ и $t\geqslant 0$, $s\geqslant 0$, при этом $X_0\neq X'_0$.
И пусть известно, что в некоторый момент времени $\tau$ эти процессы совпали: $X_\tau=X'_\tau=\mathfrak{X}$. 
Тогда, в соответствии с переходной функцией и Марковским свойством, для всех $t=\tau+s>\tau$ и $A\in \sigma(\XX)$ выполняется: $\PP\{X_t\in A\}=\PP\{X_t\in A|X_\tau=\mathfrak{X}\}=\PP\{X_t'\in A|X'_\tau=\mathfrak{X}\}=\PP\{X'_t\in A\}$.

Если известно распределение (или оценка распределения) момента склеивания $\tau$ (эта величина, вообще говоря, зависит от начальных условий процессов $X_t$ и $X_t'$), то выписывается {\it основное неравенство склеивания}:
\begin{multline}\label{3}
|\PP\{X_t\in S\}-\PP\{X_t'\in S\}|=
\\ \\
=|\PP\{X_t\in S\}-\PP\{X_t'\in S\}|\times (\1(\tau> t)+ \1(\tau\leqslant t))\leqslant \PP\{\tau> t\}    \end{multline}
(здесь $\1(\cdot)$ обозначает индикатор события).

Процессы $\left(X_t,\,t\geqslant 0\right)$ и $\left(X_t',\,t\geqslant 0\right)$ имеют распределения \linebreak $\PPP_t(A)\bd \PP\left\{X_t\in A\right\}$, и $\PPP_t'(A)\bd \PP\left\{X_t'\in A\right\}$ соответственно. 

$
\tau \left(X_0,X_0' \right)\stackrel{\rm{def}}{=\!\!\!=\!\!\!=} \inf \left\{ t>0: \,X_t=X_t' \right\} .
$
Предположим, что для некоторой положительной возрастающей функции $\varphi(t)$ вычислена или оценена сверху величина $\EE\,\varphi\left(\tau \left(X_0,X_0' \right)\right)=C \left(X_0,X_0' \right)<\infty$ .
Подставим эту функцию в последнее выражение (\ref{3}) и применим неравенство Маркова:
\begin{multline*}
\left|\PPP_t( A)-\PPP_t'( A) \right|\leqslant \mathbf P\left\{\tau \left(X_0,X_0' \right) >t \right\}= \mathbf P\left\{\varphi\left(\tau \left(X_0,X_0' \right)\right) >\varphi(t) \right\}\leqslant
\\ \\
\leqslant \left|\PPP_t( A)-\PPP_t'( A) \right|\leqslant\dd \frac{\mathbf{E}\, \varphi \left(\tau \left(X_0,X_0' \right) \right)}{\varphi(t)} . 
\end{multline*}
Если известна оценка стационарного распределения $\PPP$ процесса $X_t$ (и $X_t'$ -- это одно и то же), то последнее неравенство можно проинтегрировать по стационарной мере $\PPP$, и тогда получаем
$$
\left|\PPP_t( A)-\PPP(A) \right|\leqslant (\varphi(t))^{-1} \underbrace{\int\limits_{\mathscr{X}}\varphi \left (\tau \left (X_0,X_0' \right) \right ) \,\mathrm{d}}_{\widehat C\left (X_0\right )}\PPP \left (X_0'\right)= \frac{\widehat C\left (X_0\right )}{\varphi(t)},
$$
т.е. 
$$
\left\|\PPP_t^{X_0}-\PPP \right\|_{\mbox{ПВ}}\leqslant \frac{\widehat C\left (X_0\right )}{\varphi(t)}.
$$

Напомним
\begin{mydef}
Расстояние в метрике полной вариации между двумя вероятностными мерами $\PPP_1$ и $\PPP_2$ на $\sigma$-алгебре $\FFF$ подмножеств вероятностного пространства ${ \Omega }$ определяется как

$${\displaystyle \delta (\PPP_1,\PPP_2)=\|\PPP_1-\PPP_2\|_{\mbox{ПВ}}=\sup _{A\in {\FFF}}\left|\PPP_1(A)-\PPP_2(A)\right|.}
\eqno{\TR}$$

\end{mydef}

Метод склеивания впервые был предложен в 1938 г. молодым французским математиком немецкого происхождения В.Доблином (V.Doblin или \linebreak W.Doeblin) для обычных цепей Маркова. 
К сожалению, в отечественной учебной литературе практически нет информации об этом методе.
Он упоминается в «Дополнении» в пособии \cite{3.}; см. также \cite{9.}.

Для процессов с непрерывным временем применение метода склеивания ``напрямую'' невозможно, поскольку должны совпасть непрерывные сл.в. -- а вероятность такого совпадения равна нулю.

Т.е. в общем случае для процессов в непрерывном времени $$\PP \left\{ \tau \left (X_0,X_0' \right )<\infty \right\} <1,$$ и стандартное применение метода склеивания невозможно.
Поэтому применяется процедура конструирования {\it параллельного склеивания}.
\subsection{\it Параллельное склеивание (Successful coupling -- см. \cite{10.}}

Для оценки скорости сходимости друг к другу распределений двух {\it независимых} однородных Марковских процессов $\left(X_t,\,t\geqslant 0\right)$ и $\left(X_t',\,t\geqslant 0\right)$ с одинаковыми переходными функциями конструируется (на специально выбранном вероятностном пространстве) парный случайный процесс the $\mathscr{Z}_t= \left(\left (Z_t ,Z_t ' \right),t\geqslant0 \right)$ такой, что:

1. Для всех $t\geqslant 0$ сл.в. $X_t$ и $Z_t$ имеют одинаковое распределение, также как и пара $X_t'$ и $Z_t'$. Иначе говоря, маргинальные распределения процессов $X_t$ и $Z_t$ (или $X_t'$ и $Z_t'$) совпадают -- но это не значит, что совпадают конечномерные распределения у этих пар процессов.

2. $\mathbf{E}\,\tau \left(Z_0 ,Z_0 ' \right)<\infty\mbox{ где }\tau \left(Z_0 ,Z_0 ' \right)={\tau} (\mathscr Z_0 )\stackrel{\rm{def}}{=\!\!\!=\!\!\!=} \inf \left\{ t\leqslant 0:\,Z_t =Z_t ' \right\} 
$ -- момент склеивания процессов $Z_0$ и $Z_0 '$.

3. $Z_t =Z_t '$ для всех $t\geqslant {\tau} \left(Z_0 ,Z_0 ' \right)$.

Парный случайный процесс $\mathscr{Z}_t=\left(\left (Z_t ,Z_t ' \right),t\geqslant0 \right)$, удовлетворяющий условиям 1--3, называется is called {\it параллельным склеиванием}, а момент $\tau$ называется {\it успешной склейкой}.

.

{\it Повторимся, что конечномерные распределения процессов $\left(Z_t^{(i)},\,t\geqslant0\right)$ могут отличаться от конечномерных распределений процессов $\left(X_t^{(i)},\,t\geqslant0\right)$; более того, при конструировании параллельного склеивания процессы $Z_t$ и $Z_t'$, как правило, зависимы.}

Для всех $A\in \mathscr{B}(\mathscr{X})$ основное неравенство склеивания (\ref{3}) переписывается так:
\begin{multline}\label{4}
\left|\PPP^{X_0 }_t(A) - \PPP_t^{X_0 '}(A) \right|= \left|\PP \left\{ X_t \in A \right\} - \PP \left\{ X_t '\in A \right\} \right|= 
\\
\\
= \left|\PP \left\{ Z_t \in A \right\} - \PP \left\{ Z_t '\in A \right\} \right| \leqslant
\PP \left\{ {\tau}\left (Z_0 ,Z_0 '\right )\geqslant t \right\}
\leqslant
\\
\\
\leqslant\frac{\mathbf {E} \,\varphi \left({\tau} \left(Z_0 ,Z_0 ' \right) \right)}{\varphi(t)}\leqslant \frac{C \left(Z_0 ,Z_0 ' \right) }{\varphi(t)},
\end{multline}
где значение $C \left(Z_0 ,Z_0 ' \right)$ вычисляется по начальным условиям процесса $\mathscr{Z}_t$.

Но $Z_0=X_0$ и $Z_0'=X_0'$, и поэтому правая часть неравенства (\ref{4}) зависит только от $(X_0,X_0')$.

Поэтому если $\PPP_t^{X_0 }\Longrightarrow \PPP $ для всех начальных состояний $X_0$, можно использовать интегрирование по стационарной мере $\PPP $ как и выше. 

{\it Хотя в ряде случаев это интегрирование вызывает некоторые трудности.}

\section{Конструирование параллельного склеивания для перескоков процесса восстановления}

В этом разделе приводится базовая схема конструирования параллельного склеивания, на основе которой можно строить оценки скорости сходимости распределений различных регенерирующих процессов.

Итак, мы конструируем параллельное склеивание
$$(\mathscr Z_t,\,t\geqslant0)=\left(\left(Z_t , Z_t '\right),\,t\geqslant0\right)$$
для перескоков вложенных процессов восстановления $\left(B_t ,\,t\geqslant0\right)$ и $\left( B_t ',\,t\geqslant0\right)$ для исследуемых регенерирующих процессов $X_t$, $X_t'$ с различными начальными состояниями.
Это значит, что процессы $B_t $, и $B_t '$ имеют разные начальные состояния.

Т.е. соответствующие процессы восстановления начинаются {\it до} момента $t=0$, и в момент $t=0$ значения перескоков вложенных процессов восстановления различны: $B_0 ={b}$, $B_0 '=b'$; что можно обозначить как $B_t =B_t^{b}$, $B_t '=B_t^{b'}$.

Распределения этих процессов в момент $t$ обозначим соответственно $\PPP^{b}_t$ и $\PPP^{b'}_t$.

Если мы оценим сл.в. $\tau(b,{b'})=\tau(b,{b'})\stackrel {\text{\rm def}}{=\!\!\!=} \inf\left\{ t>0:\, Z_t = Z_t '\right\}$, то можно получить оценку  
$$\left\|\PPP^{b}_t(A) - \PPP^{b'}(A) \right\|_{\mbox{ПВ}}\leqslant \PP \left\{ \tau({b},{b'})>t\right\}\leqslant \frac{\EE\, \varphi( \tau({b},{b'}))}{\varphi(t)}.
$$
\subsection{\it Некоторые сведения}

В ходе конструирования параллельного склеивания надо уметь конструировать случайные величины (на некотором вероятностном пространстве).
\begin{myrem}
Напомним, как конструировать сл.в. по её ф.р. $F(s)$.

Определим обратную функцию к ф.р. $F(s)$ как 
$$
F^{-1}(y)\bd\inf\{x:\,F(x)\geqslant y\}.
$$
Если $\UU$ -- равномерно распределённая на $[0,1)$ сл.в., то $\xi\bd F^{-1}(\UU)$ имеет ф.р. $F(s)$. \TR
\end{myrem}

Ранее говорилось, что  $B_t$ -- это перескок процесса восстановления $N_t$ с ф.р. времени восстановления $F(s)$; 
$N_t\stackrel {\text{\rm def}}{=\!\!\!=} \sum\limits_{i=0}^\infty \mathbf{1} (S_n<t)$ где $S_n\stackrel {\text{\rm def}}{=\!\!\!=} \sum\limits_{n=0}^n \xi_i$, 
$\mathbf {E} \,\xi_i<\infty$ и $\xi_i$, $i\in \mathbb Z_+$ независимы в совокупности; 
$\PP \left\{ \xi_i\leqslant s\right\} =F(s)$, $i\in\mathbb N$.
$B_t\stackrel {\text{\rm def}}{=\!\!\!=} (t-S_{N_t})$ -- перескок процесса восстановления $N_t$.

Если $B_0=b$ при $t=0$, значит, последнее восстановление процесса $N_t$ было в момент $t=-b$, и остаточное время периода восстановления (недоскок) $W_0$ имеет ф.р. 
\begin{multline}\label{5}
\PP\{W_0\leqslant s|B_0=b\}= \PP\{\xi_0-b\leqslant s|\xi_0>b\}=\PP \left\{\xi_0\in(b,b+s]| \xi_0>b\right\} =
\\ \\
=\dd\frac{F(s+b)-F(b)}{1-F(b)} \stackrel {\text{\rm def}}{=\!\!\!=}F_b(s) .
\end{multline}

Также напомним, что 
\begin{equation}\label{6}
\lim\limits_{t\to \infty}\PP \left\{ B_t\geqslant s\right\} =\mu^{-1}\dd\int\limits_{s}^\infty (1-F(u))\,\mathrm{d}\, u\stackrel {\text{\rm def}}{=\!\!\!=} 1-\widetilde F(s),    
\end{equation}
где
$\mu\stackrel {\text{\rm def}}{=\!\!\!=} \mathbf {E} \,\xi_1$ (см. \cite{8.}).

\subsection{\it Первый шаг: конструирование процессов восстановления}

Везде в дальнейшем $\UU_i$, $\UU_i'$, $\UU_i''$, $\UU_i'''$,\ldots -- независмые равномерно распределённые на $[0,1)$ сл.в..

1. Конструирование процесса $(B _t,\,t\geqslant 0)$с начальным состоянием $B_0 =b$.
Обозначим $t_1\bd F_{b}^{-1}(\UU_1)$, $t_k\bd t_{k-1}+F^{-1}(\UU_k)$, $k>1$.

\begin{figure}[h]
    \centering
\begin{picture}(350,50)
\put(0,20){\vector(1,0){350}}
\thicklines
{\qbezier(0,33)(40,34)(60,20)}
{\qbezier(60,20)(80,40)(100,20)}
\qbezier(100,20)(140,40)(180,20)
\qbezier(180,20)(200,40)(220,20)
\qbezier(220,20)(250,40)(280,20)
\qbezier(280,20)(290,40)(300,20)
\qbezier(300,20)(320,31)(330,30)
\thinlines
\put(5,37){\small $\xi_0\sim F_b$}
\put(65,37){\small $\xi_1\sim F$}
\put(125,37){\small $\xi_2\sim F$}
\put(185,37){\small $\xi_3\sim F$}
\put(235,37){\small $\xi_4\sim F$}
\put(210,11){\small $t_3$}
\put(58,11){\small $t_0$}
\put(98,11){\small $t_1$}
\put(178,11){\small $t_2$}
\put(278,11){\small $t_4$}
\put(298,11){\small $t_5$}
\put(0,11){\small $0$}
\end{picture}
    \caption{Процесс ${{X}_{t}}$ (соответствующий ${{B}_{t}}$).}
    \label{fig3}
\end{figure}

2. Аналогично для $B_t '$ с начальным состоянием $ B_0 '=b'$ положим $t_1'\bd F_{b'}^{-1}(\UU_1')$, $t_k'\bd t_{k-1}'+F^{-1}(\UU_k')$, $k>1$.

\begin{figure}[h]
    \centering
\begin{picture}(350,50)
\put(0,20){\vector(1,0){350}}
\thicklines
{\qbezier(0,33)(30,34)(40,20)}
{\qbezier(40,20)(60,40)(80,20)}
\qbezier(80,20)(120,40)(160,20)
\qbezier(160,20)(200,40)(240,20)
\qbezier(240,20)(260,40)(280,20)
\qbezier(280,20)(320,31)(330,30)
\thinlines
\put(0,37){\small $ \xi_0'\sim F_{b'}$}
\put(45,37){\small $\xi_1'\sim F$}
\put(105,37){\small $\xi_2'\sim F$}
\put(185,37){\small $\xi_3'\sim F$}
\put(255,37){\small $\xi_4'\sim F$}
\put(238,11){\small $ t_3'$}
\put(38,11){\small $ t_0'$}
\put(78,11){\small $ t_1'$}
\put(158,11){\small $ t_2'$}
\put(278,11){\small $ t_4'$}
\put(0,11){\small $0$}
\end{picture}
    \caption{Процесс ${{{X}'}_{t}}$ (соответствующий ${{{B}'}_{t}}$).}
    \label{fig4}
\end{figure}

Понятно, что вероятность склеивания (т.е. совпадения двух точек $t_i$ и $t_j$) за конечное время равна нулю.

Более того, построенные процессы независимы, а в параллельном склеивании они будут зависеть.

Для конструирования параллельного склеивания (т.е. модификации предложенной выше конструкции) нам понадобится Основная Лемма Склеивания (см., например, \cite{12.}).

\begin{mylem} [Основная Лемма Склеивания]
Пусть $f_i(s)$ -- плотности распределения (п.р.) сл.в. $\theta_i$ ($i=1,2$).
и $$\intl_{-\infty}^ \infty \min(f_1(s),f_2(s))\ud s=\varkappa>0.$$

Тогда на некотором вероятностном пространстве существует две сл.в. $\vartheta_i$ такие, что $\vartheta_i\bD \theta_i$, и $\PP\{\vartheta_1=\vartheta_2\}\geqslant \varkappa$. \hfill \ensuremath{\triangleright}
\end{mylem}
Величина $\dd\intl_{-\infty}^ \infty \min(f_1(s),f_2(s))\ud s=\varkappa$ называется {\it общей частью} распределений $\theta_i$.
\begin{proof}
Пусть 
$$\varphi(s)\stackrel {\text{\rm def}}{=\!\!\!=} \min\left(f_1(s), f_2(s)\right); \qquad \dd\int\limits_0^\infty \varphi(s)\,\mathrm{d}\, s\stackrel {\text{\rm def}}{=\!\!\!=} \varkappa=\Phi(+\infty)>0,$$
где $\Phi(s)\stackrel {\text{\rm def}}{=\!\!\!=} \dd\int\limits_0^s \varphi(u)\,\mathrm{d}\, u$.

Обозначим
$$\Psi(s)\stackrel {\text{\rm def}}{=\!\!\!=} F_1(s)-\Phi(s), \;\;  \Psi_c(s)\stackrel {\text{\rm def}}{=\!\!\!=} F_2 (s)-\Phi(s); \;\; \Psi_1(+\infty)=\Psi_2(+\infty)=1-\varkappa,
$$
где $F_i(s)$ -- ф.р. сл.в. $\theta_i$.

Положим для независимых равномерно распределённых на $[0,1)$ сл.в. $\UU$, $\UU'$, $\UU''$ 
$$
\xi_1(\UU,\UU',\UU'')\stackrel {\text{\rm def}}{=\!\!\!=} \mathbf{1} (\UU<\varkappa)\Phi^{-1}(\varkappa \UU')+\mathbf{1} (\UU\geqslant \varkappa) \Psi^{-1}((1-\varkappa) \UU'');$$
$$\xi_2(\UU,\UU',\UU'')\stackrel {\text{\rm def}}{=\!\!\!=} \mathbf{1} (\UU<\varkappa)\Phi^{-1}(\varkappa \UU')+\mathbf{1} (\UU\geqslant \varkappa) \Psi_c^{-1}((1-\varkappa) \UU'').
$$

Несложно заметить, что
$$\PP\Big\{\xi_1(\UU,\UU',\UU'')\leqslant s\Big\}=F_1(s);\quad \PP\Big\{\xi_2(\UU,\UU',\UU'')\leqslant s\Big\}= F_2 (s);
$$
и $\PP\left\{\xi_1(\UU,\UU',\UU'')=\xi_2(\UU,\UU',\UU'')\right\}=\varkappa$.
\end{proof}
\begin{myrem}
Фактически 
$$
\xi_i=
\begin{cases}
\widetilde \xi & \mbox{с вероятностью } \varkappa;
\\ \\
\widehat \xi_i & \mbox{с вероятностью } 1-\varkappa,
\end{cases}
$$
где независимые сл.в. $\widetilde \xi$ и $\widehat \xi_i$ имеют ф.р. $\dd\frac{\Phi(s)}{\varkappa}$ и $\dd\frac{\Psi(s)}{1-\varkappa}$. \TR
\end{myrem}

\subsection{\it Второй шаг: модификация конструирования пары $(B_t,B_t '$) в параллельное склеивание}

Вернёмся к нашим {\it независимым} процессам $B_t $ и $B_t '$.
У нас $B_0 =b$, $B_0 '=b'$, это -- перескоки процессов восстановления.

Т.е. времена $t_1$ и $t_1'$ -- остаточные времена периодов регенерации с ф.р. $F_{b}(s)$ и $F_{b'}(s)$ соответственно (см. (\ref{5})). 
Первый этап конструирования -- это построение {\it независимых} сл.в. $t_1$ и $t_1'$.

После момента $T_1\bd \max\left\{t_1 ,t_1 '\right\}$ оба процесса начинают удовлетворять условиям Теоремы Лордена (очевидно, $T_1\leqslant t_1 +t_1 '\bd \theta_0$).

На рисунке \ref{fig5} $T_1= t_1 '$. После момента $T_1$ к процессу $B_t'$ можно применять неравенство Лордена, т.е. $\EE\,B_{t '_1} \leqslant\Xi$.

\begin{figure}[b]
\centering
\begin{picture}(450,100)
\put(0,20){\vector(1,0){450}}
\put(0,18){\line(0,1){4}}
\thicklines
{\qbezier(0,35)(50,35)(100,20)}
{\qbezier(100,20)(140,40)(180,20)}
\qbezier(180,20)(240,40)(300,20)
\qbezier(300,20)(340,40)(380,20)
\put(400,40){\ldots \ldots \ldots}
\put(95,5){$t _1$}
\put(175,5){$t _2$}
\put(298,5){$t _3$}
\put(378,5){$t _4$}
\put(40,40){\small $\xi _1=t _1$}
\put(137,40){\small $\xi _2$}
\put(220,40){\small $\xi _3$}
\put(330,40){\small $\xi _4$}
\put(-5,5){$t_0=0$}
\put(243,5){$t_2 '$}
\bk
\put(0,65){\vector(1,0){450}}
\put(0,63){\line(0,1){4}}
\thicklines
{\qbezier(0,80)(80,80)(150,65)}
{\qbezier(150,65)(200,85)(243,65)}
\qbezier(243,65)(269,85)(300,65)
\qbezier(300,65)(340,85)(380,65)
\put(150,75){$t '_1$}
\put(243,75){$t '_2$}
\put(300,75){$t '_3$}
\put(378,75){$t '_4$}
\put(40,85){\small$ \xi '_1=t '_1$}
\put(200,85){\small$ \xi '_2$}
\put(269,85){\small $\xi '_3$}
\put(330,85){\small$ \xi '_4$}
\put(-5,50){$t_0=0$}
\gr
\multiput(150, 5)(0,5){15}{\line(0,1){3}}
\bk
\put(150,5){$t '_1$}
\rr
%
\qbezier(180,20)(210,0)(243,20)
\put(120,-1){\small $B_{t '_1} $}
\rr
%
\qbezier(100,20)(125,0)(150,20)
\put(210,-1){\small $B_{t '_2} $}

\gr
\multiput(243, 5)(0,5){15}{\line(0,1){3}}
\bk
\thinlines
\bb
\multiput(300, 5)(0,5){15}{\line(0,1){3}}
\put(305,55){\small\it склейка}
\put(310,50){\vector(-1,-1){10}}
\bk
\end{picture}
 \caption{На этом рисунке в точке ${{{t}'}_{2}}$ происходит склеивание.}
    \label{fig5}
\end{figure}

По неравенству Маркова, для любого $\Theta>\Xi\geqslant \EE\,B_{t '_1}'$ верно неравенство $\dd\PP\left\{B_{t '_1} '<\Theta\right\}\geqslant 1-\frac{\Xi}{\Theta}=p_0$ и в случае, если произойдёт событие $B_{t '_1} <\Theta$, можно оценить величину
$$\dd\supl_{u\in (0,\Theta)}\dd\intl_{0}^{\infty} \min(\lambda(s+u),\lambda(s))\ud s\geqslant \supl_{u\in (0,\Theta)}\intl_{0}^{\infty} \min(\varphi(s+u),\varphi(s))\ud s=\kappa>0,$$
т.к. почти всюду $\varphi(s)>0$ (поскольку предполагается, что плотность распределения сл.в. $\xi_i$ п.в. положительна).

Теперь, используя Основную Лемму Склеивания, можно продолжить остаточное время периода восстановления процесса $B_t$ и следующий период восстановления $B '_t$ таким образом, что с вероятностью $\varkappa\bd p_0\kappa$, процессы $B _t$ и $B '_t$ совпадут в момент $t_2 '\bd \theta_1$.

(Заметим, что здесь использовалась некоторая величина $\Theta>\Xi$; выбор этой величины влияет на значение $\varkappa$ -- так что можно искать значение $\Theta$, при котором значение $\varkappa$ максимально.) 

Если в момент $t_2 '=\theta_1$ не произошло совпадения процессов (с вероятностью $\leqslant (1-\varkappa)$), то в следующий момент восстановления $B '_t$ мы повторяем описанную выше конструкцию с перескоком $B _t$ процесса восстановления в момент $t_2 '=\theta_2$ таким образом, что процессы $B _t$ и $B '_t$ совпадут с вероятностью $p_0\varkappa$ в момент $t_3 '=\theta_3$, и т.д..
Таким образом, при условии, что в момент $t_n '=\theta_{n-1}$ произойдёт совпадение процессов $B _t$ и $B '_t$,
$$\tau\left(B_0 ,B_0^{(t)}\right)\leqslant T_1+\dd\suml_{i=2}^n \xi_i '$$ 
с вероятностью $\varkappa(1-\varkappa)^{n-2}$ (здесь $T_1\leqslant t_1 +t_1 '$).

Иначе говоря, момент склеивания $\tau\left(B_0 ,B_0^{(t)}\right)$ оценивается сверху {\it условной} геометрической суммой периодов восстановления, включая сумму первый (возможно, неполных) периодов восстановления $B_t$ и $B '_t$.

Если существует конечный момент $\EE\,(\xi)^{\ell+1}=C_\ell\left(B_0 ,B_0'\right)$, то можно получить оценку сверху для 
$$
\EE\,\left(\tau\left(B_0 ,B_0^{(t)}\right)\right)^\ell=C_\ell\left(B_0 ,B_0 '\right)=C_\ell\left(b,b'\right).
$$
Для этого можно использовать хорошо известное \\
{\bf Неравенство Йенсена}

Для действительнозначной выпуклой (вниз) функции $\varphi$, чисел $ x_{1},x_{2},\ldots ,x_{n}$ из её непрерывной области определения, и положительных величин $a_{i}$, неравенство Йенсена можно сформулировать так:
$$
\varphi \left({\frac {\sum a_{i}x_{i}}{\sum a_{i}}}\right)\leq {\frac {\sum a_{i}\varphi (x_{i})}{\sum a_{i}}}
$$
и, следовательно, 
$$ 
\varphi \left({\frac {\sum x_{i}}{n}}\right)\leq {\frac {\sum \varphi (x_{i})}{n}}.
$$

Если существует $\EE\,\xi^{\ell+1}<\infty$, оценим
\begin{multline*}
\EE\,\left(T_1+\suml_{i=1}^{\nu} \xi_i\right)^\ell=\EE\,\left(\frac{(\nu+1)\times\left(T_1+\suml_{i=1}^{\nu} \xi_i\right)}{\nu+1}\right)^\ell\leqslant 
\EE\,(\nu+1)^\ell\frac{\left(T_1\right)^\ell+\suml_{i=1}^\nu \xi_i^\ell}{\nu+1}=
\\
=\EE(\nu+1)^{\ell-1}\times\EE\,T_1^\ell +\suml_{i=1}^\infty \left((i+1)\times\suml_{j=1}^{i-1}\EE\,(\xi_j^\ell|\EEE_i)\PP(\EEE_i)+\EE\,(\xi_i^\ell|\EEE_i)\PP(\EEE_i) \right),
\end{multline*}
где $\EEE_i=\bigcup\limits_{j=1}^{i-1}\overline \SSS_j\cup\SSS_i$, а $\SSS_i$ -- событие $\{$момент $\theta_i$ -- это момент склеивания $\tau\}$.

Заметим, что $\PP(\SSS_i)= (1-\varkappa)^{i-1}(1-\varkappa)\leqslant (1-\varkappa)^{i-1}$, а также используем очевидное неравенство 
\begin{equation}\label{7}
\EE\,(\xi|A)\PP(A)\leqslant \EE\,(\xi).    
\end{equation}

Поэтому $\EE\,(\xi_j^\ell|\EEE_i)\PP(\EEE_i)\leqslant\EE\,(\xi_j^\ell|\EEE_i)\prod\limits_{j=1}^{i-1}(1-\PP(\SSS_j))\leqslant \EE\, \xi_j^\ell (1-\varkappa)^{i-2}$; $\EE\,(\xi_i^\ell|\EEE_i)\PP(\EEE_i)\le \EE\,\xi_i^\ell (1-\varkappa)^{i-1}\le (1-\varkappa)^{i-2}$.

Таким образом, $$\dd\EE\,\left(T_1+\suml_{i=1}^{\nu} \xi_i\right)^\ell\le \EE(\nu+1)^{\ell-1}\times\EE\,T_1^\ell+\EE\,\xi_i^\ell \suml_{i=1}^\infty (i+1)^2(1-\varkappa)^{i-1}=\Upsilon(b,b',\Theta, F(\cdot)),$$ 
эта величина $\Upsilon(\ell,b,b',\Theta, F(\cdot))$ может быть вычислена или оценена.

Учитывая, что предельное (стационарное) распределение перескока процесса восстановления известно (\ref{5}), и только величина $T_1$ зависит от $b$ и $b'$, интегрирование величины $\Upsilon(\ell,b,b',\Theta, F(\cdot))$ не представляет больших трудностей.
Таким образом, верна 
\begin{theorem}
Если конечна величина $\EE\,\xi^{\ell+1}$, то можно вычислить величину $\dd\widetilde \Upsilon(\ell,b,b',\Theta, F(\cdot))=\intl_0^\infty \Upsilon(\ell,b,\Theta, F(\cdot)) \ud \widetilde F(b')$ такую, что
$$
\left\|\PPP^{b}_t - \widetilde\PPP \right\|_{\mbox{ПВ}}\leqslant \frac{ \widetilde \Upsilon(\ell,b,\Theta, F(\cdot))}{t^\ell}.
\eqno{\triangleright}$$
\end{theorem}

Оценка экспоненциальных моментов сл.в. $\tau(b,b')$ существенно сложнее, т.к. неравенство (\ref{7}) не даёт возможности вычислить оценку $\EE\exp(\beta \tau)$; здесь для получения оценки можно использовать свойства конкретного распределения $F(s)$ или более тонкие оценки при конструировании параллельного склеивания, на чём мы здесь не будем останавливаться. 
Тем не менее, и если $\EE\,\exp(\alpha\xi)$, то можно оценить 
$$
\EE\,\exp\left(\beta \times \left(T_1+\suml_{i=1}^{\nu} \xi_i\right)\right)
$$
для некоторого $\beta\in(0, \alpha)$, и поэтому верна
\begin{theorem}
Если $\EE\,\exp(\alpha \Xi)<\infty$, то для некоторого $\beta\in (0,\alpha)$ можно вычислить постоянную $K_\beta(b,b')$ такую, что
$\EE\,\exp(\beta \tau(b,b'))\leqslant K_\beta(b,b')$,
и
$$
\left\|\PPP^{b}_t - \widetilde\PPP \right\|_{\mbox{ПВ}}\leqslant \widetilde K_\beta(b)\exp(-\beta t),
$$
где $\widetilde K_\beta(b)=\dd\intl_0^\infty \exp(K_\beta(b,b')) \ud \widetilde F(b')$.\TR
\end{theorem}

\section{Заключение}

Использование предложенных способов вычисления строгих оценок сверху для регенерирующих процессов может быть использован при анализе поведения различных сложных СМО, СеМО и систем надёжности.
Оценив время $\mathcal T$ достижения приемлемого уровня близости распределения исследуемой модели, исследователь может с помощью имитационного моделирования проанализировать поведение системы до достижения времени $\mathcal T$ и использовать эти данные при прогнозировании поведения аналогичных систем.

Для случая, когда исследуемая модель включает несколько параллельных технологических процессов, может применяться естественное обобщение Основной Леммы Склеивания:
\begin{mylem}[Обобщение Основной Леммы Склеивания]
Пусть $f_i(s)$ -- плотности распределения сл.в. $\theta_i$ ($i=1,\ldots,n$).
и пусть $$\dd\intl_{-\infty}^ \infty \min\limits_{i=1,\ldots,n}(f_i(s))\ud s=\varkappa>0.$$

Тогда на некотором вероятностном пространстве существует $n$ сл.в. $\vartheta_i$ ($i=1,\ldots,n$), таких, что $\vartheta_i\bD \theta_i$ для всех $i=1,\ldots,n$, и $$\PP\{\vartheta_1=\vartheta_2=\ldots=\vartheta_n\}\geqslant \varkappa. \hfill \eqno{\triangleright}$$
\end{mylem}
Доказательство этой Леммы конструктивно и подобно доказательству Основной Леммы Склеивания  (см., например, \cite{13.}).

Используя Обобщение Основной Леммы склеивания и предложенный в \cite{13.} подход к конструированию многомерных марковски модулированных процессов, можно получать оценки момента склеивания (склейки) сконструированных по предложенной выше схеме ``параллельных’’ (т.е. имеющих те же маргинальные распределения) случайных процессов для широкого круга задач
ТМО, СеМО и теории надёжности.

Умение оценивать сверху, а не только указывать тип скорости сходимости (степенной, экспоненциальный) важно для практического использования различных стохастических моделей.
Использование таких оценок позволяет узнать, когда распределение используемой системы станет достаточно близко к предельному (стационарному) распределению. 
До этого момента можно оценивать поведение используемой системы с помощью имитационного моделирования, что не всегда бывает экономически оправдано.

Поэтому так важно уметь оценивать {\it гарантированное} время до приближения распределения используемой системы к стационарному режиму.

\section*{Благодарности }
Авторы благодарят анонимного рецензента за ценные замечания, способствовавшие улучшению текста.
Работа частично поддержана РФФИ (проекто № 20-01-00575А)

\end{document}